\documentclass[11pt,a4paper]{article}

\usepackage{amsfonts,amsgen,amsmath,amstext,amsbsy,amsopn,amsthm,amssymb
}
\usepackage{authblk}
\usepackage{color}
\usepackage{verbatim}
\usepackage{multirow}
\usepackage[pdfencoding=auto, psdextra]{hyperref}
\usepackage[noabbrev,capitalise]{cleveref}
\usepackage[title]{appendix}
\usepackage{tikz}
\usepackage{arydshln}

\newtheorem{dfn}{Definition}[section]

\newtheorem{thm}[dfn]{Theorem}
\newtheorem{lem}[dfn]{Lemma}
\newtheorem{prop}[dfn]{Proposition}

\newtheorem{clm}[dfn]{Claim}

\numberwithin{equation}{section}

\newcommand{\mZ}{\mathbb{Z}}
\newcommand{\eps}{\varepsilon}
\newcommand{\cok}{\operatorname{cok}}
\newcommand{\im}{\operatorname{im}}
\newcommand{\vvert}{\,\vert\,}

\newcommand{\rb}{\tilde{b}}
\newcommand{\rH}{\tilde{H}}

\def\Xv#1{X\cup\{#1\}\vvert Y}
\def\vY#1{X\vvert Y\cup\{#1\}}
\def\XvY#1#2{X\cup\{#1\}\vvert Y\cup\{#2\}}

\def\qed{\ifhmode\unskip\nobreak\hfill$\Box$\bigskip\fi \ifmmode\eqno{Box}\fi}

\begin{document}

\title{\bf\Large The Betti Number of the Independence Complex of Ternary Graphs}

\author[1]{Wentao Zhang\thanks{Email: {\tt wtzhang20@fudan.edu.cn}}}

\author[1] {Hehui Wu\thanks{Supported in part by National Natural Science Foundation of China grant 11931006, National Key Research and Development Program of China (Grant No. 2020YFA0713200), and the Shanghai Dawn Scholar Program grant 19SG01.Email:{\tt hhwu@fudan.edu.cn.}}}

\affil[1]{Shanghai Center for Mathematical Sciences\\
Fudan University\\Shanghai, China}

\maketitle

\begin{abstract}
Given a graph $G$, the \textit{independence complex} $I(G)$ is the simplicial complex whose faces are the independent sets of $V(G)$. Let $\rb_i$ denote the $i$-th reduced Betti number of $I(G)$, and let $b(G)$ denote the sum of $\rb_i(G)$'s. A graph is ternary if it does not contain induced cycles with length divisible by three. G. Kalai and K. Meshulam conjectured that $b(G)\le 1$ whenever $G$ is ternary. We prove this conjecture. This extends a recent results proved by Chudnovsky, Scott, Seymour and Spirkl that for any ternary graph $G$, the number of independent sets with even cardinality and the independent sets with odd cardinality differ by at most 1.
\end{abstract}

\section{Introduction}

A graph is \textit{ternary} if it has no induced cycle of length divisible by three. A ternary graph is also called a Trinity Graph by others \cite{BCT, Ka14}. Given a graph $G$, let  $f_G$ be the sum of $(-1)^{\vert A\vert}$ over all independent sets $A$. Recently, Chudnovsky, Scott, Seymour and Spirkl~\cite{CSSS20} proved a intriguing conjecture on the independent sets (or stable sets) of ternary graphs proposed by G. Kalai and R. Melshulam (see \cite{Ka14}) in the late 1990's.

\begin{thm}\label{thm:EC}
	If $G$ is a graph with no induced cycle of length divisible by three,  then $\vert f_G\vert\leq1$.
\end{thm}

A stronger version of the conjecture of Kalai and Meshulam concerns the Betti number of the independence complex of a ternary graphs, which build a connection between algebraic topology and graph theory.

The \textit{independence complex} $I(G)$ of a graph $G$ is the simplicial complex whose faces are the independent sets of $V(G)$. $\rH_i(I(G))$ is the $i$-th reduced homology group of $I(G)$, and $\rb_i(I(G))=dim \rH_i(I(G))$ is the $i$-th reduced Betti number of $I(G)$. Note that the Betti number $b_i$ of a simplicial complex equals the reduced Betti number, except only for the 0-th Betti number, which is one more than $\rb_0$. Specially, when $G$ is a null graph (with no vertex), we let $b_0(I(G))=0$ and $\rb_0(I(G))=-1$.

Let $b(G)$ denote the sum of $\rb_i(I(G))$'s. For a simplicial complex, the Euler characteristic can be defined as $\sum_i (-1)^{i} b_i$, which is $1+\sum_i (-1)^i\rb_i$. From a basic theorem in homology theory, we know that the Euler characteristic of $I(G)$ also equals $\sum (-1)^{|A|-1}$, over all the non-empty independent sets in $G$ (see \cite{Ha09}). It immediately follows that $f_G=\sum_{i=0}^{\infty}(-1)^{i+1}\rb_i(G)$, and so $\vert f_G\vert\leq b(G)$. 
 
 Note that $b(G)\ge |f_G|=2$ when $G$ is an cycle of length divisible by 3. A question was asked by Kalai and Meshulam (see \cite{Ka14}) on the betti number of graphs without induced cycle of length divisible by 3. The purpose of the paper is to prove this conjecture of Kalai-Meshulam (see \cite{Ka14}), which is a stronger version of Theorem~\ref{thm:EC}.

\begin{thm}\label{thm:main}
	If $G$ is a graph with no induced cycle of length divisible by three,  then $b(G)\leq1$.
\end{thm}

Analogously, a \textit{clique complex} of a graph $G$ is the simplicail complex whose faces are the cliques of $G$. In an abstract simplicial complex, a set $S$ of vertices that is not itself a face of the complex, but such that each pair of vertices in S belongs to some face in the complex, is called an \textit{empty simplex}. A \textit{flag complex} is an abstract simplicial complex that has no empty simplex. As any flag complex is the clique complex of its 1-skeleton, and  the clique complex of a graph $G$ is the independence complex of the complement of $G$, the above theorem give a full characterization of minimal flag complex with total Betti number 2.

If we further forbid any $C_{3k}$ as a subgraph instead of a induced subgraph, the following results has been claimed by A. Engstrom~\cite{En14}, which extend a result of Gauthier~\cite{Ga17} on $f_G$ of such graphs:

\begin{thm}
If $G$ is a graph without cycles of length divisible by three, then $I(G)$ is contractible or homotopy equivalent to a sphere.
\end{thm}

In the same paper, Engstrom also asked whether for any ternary graph, $I(G)$ is contractible or homotopy equivalent to a sphere. Very recently, based on our proof, J. Kim~\cite{Kim21} claimed that Engstrom's conjecture is also true.

There are some other conjectures asked simultaneously by Kalai and Meshulam since 1990's (see \cite{Ka14}), relating chromatic numbers, Euler Characteristc or Betti number of the independence complex, and ternary graphs. Some of them have been answered recently. See the paper of M. Bonamy, P. Charbit, S. Thomass\'{e}\cite{BCT},  the paper of A. Scott and P. Seymour\cite{SS19}, and some other papers (\cite{CS94, En20}).

Our proof is inspired by the proof of Theorem~\ref{thm:EC} by Chudnovsky, Scott, Seymour and Spirkl\cite{CSSS20}. The proof can be shorter if we use their results directly. But here we prefer to give an full and independent proof, as the Betti Number will give you more detail on the induction process than the Euler Characteristic, and the proof are smoother and shorter than their original paper after the system is set up.

Given a graph $G$ and vertex sets $X, Y$ let $f_G(X,Y)$ be the sum of $(-1)^{|A|}$, with $A$ goes over all the independent sets that includes $X$ and are disjoint from $Y$. Their proof is based on the recursive formula of $f_G$: $f_G(X,Y)=f_G(X\cup\{v\},Y)+f_G(X, Y\cup\{v\})$ for every $v\in V(G)$. To recursively calculate $\rb_i(I(G))$, we will instead use a formula from the Mayer-Vietoris Sequence, which is a powerful tool in calculation of homology group.

\section{Mayer-Vietoris Sequence}\label{section:Pre}

As some graph theoryists may not be very family with homology theory, we first introduce some prerequisites from homology theory. This part we refer to a paper of Delfinado and Edelsbrunner~\cite{DE95}. 

An abstract simplicial complex $K$ is a family of sets that is closed under taking subsets. Each element in the set is called a vertex and each finite set in $K$ is called a face. An $n$-face of $K$ is a face with size $n+1$. Each $n$-face can be oriented with a linear order of its vertices, denoted by $[v_0,\dots v_n]$. The chain group $C_n(K)$ is the free abelian group generated by oriented $n$-faces of $K$, and the boundary map $\partial_n:C_n(K)\to C_{n-1}(K)$ is defined by
\[
\partial_n[v_0,\dots,v_n] = \sum\limits^n_{j=0}(-1)^j[v_0,\dots,\hat{v}_j,\dots,v_n]\]
where $\hat{v}_j$ means $\hat{v}_j$ is omitted.

The \textit{reduced homology groups} $\rH_i(K)$ is the homology groups of the augmented chain complex
\[\cdots\longrightarrow C_2(K)\stackrel{\partial_2}\longrightarrow C_1(K)\stackrel{\partial_1}\longrightarrow C_0(K)\stackrel{\eps}\longrightarrow\mZ\longrightarrow0\]
where $\eps:C_0(K)\to\mZ$ is the augmentation map defined by $\eps(v)=1$ for each vertex $v$ of $K$. We will have $\partial_{i}\circ\partial_{i+1}=0$ and $\eps\circ\partial_1=0$. The $i$-th reduced homology group $\rH_i(K)$ of $K$ is the quotient group $ker(\partial_{i})/im(\partial_{i+1})$ for positive $i$, and $\rH_0(K)=ker(\eps)/im(\partial_1)$. The $i$-th reduced Betti number $\tilde{\beta}_i(K)$ of $K$ is the dimension of $\rH_i(K)$.

Let $K'$ and $K''$ be subcomplexes such that $K=K'\cup K''$ and let $L=K'\cap K''$. A chain complex is exact if $im(\partial_{i+1})=ker(\partial_{i})$ for all $i$. There is an exact sequence of reduced homology groups called the Mayer-Vioteris sequence~\cite{Mun84}
\[
\begin{split}
\cdots\longrightarrow\tilde{H}_i(L)\stackrel{\lambda_i}\longrightarrow\tilde{H}_i(K')\oplus\tilde{H}_i(K'')\longrightarrow\tilde{H}_i(K)\longrightarrow\tilde{H}_{i-1}(L)\stackrel{\lambda_{i-1}}\longrightarrow\tilde{H}_{i-1}(K')\oplus\tilde{H}_{i-1}(K'')&\longrightarrow\cdots\\\cdots\longrightarrow\tilde{H}_0(&K)\longrightarrow0
\end{split}
\]
It follows that 
\[0\longrightarrow\cok\lambda_i\longrightarrow\tilde{H}_i(K)\longrightarrow\ker\lambda_{i-1}\longrightarrow0\]
is a short exact sequence of abelian groups. Let $N_{i}=\ker\lambda_{i}$ and $\beta(N_{i})$ be its dimension. We have (see \cite{AH35}):
\begin{equation}\label{recursion}
\begin{split}
\beta_i(K)&=\beta(\cok\lambda_i)+\beta(\ker\lambda_{i-1})\\
&=\beta(\tilde{H}_i(K')\oplus\tilde{H}_i(K'')/\im\lambda_i)+\beta(N_{i-1})\\
&=\beta_i(K')+\beta_i(K'')-\beta_i(L)+\beta(N_i)+\beta(N_{i-1})
\end{split}
\end{equation}
Note that $\beta(N_i)\le \beta_i(L)$ as $N_i\subseteq \tilde{H}_i(L)$ for each $i$.

Given a graph $G$, suppose $X$ is an independent set of $G$ and $Y$ is a vertex set disjoint from $X$, let $G(X\vvert Y)$ be the subgraph induced by $V(G)-N[X]-Y$. If the elements of $X$ or $Y$ are listed we omit the braces for simplicity. For simplicity, when $G(X\vvert Y)$ is not a null graph, we write $I(X\vvert Y)$ and $b(X\vvert Y)$ for $I(G(X\vvert Y))$ and $b(G(X\vvert Y))$ when $G$ is known. Similarly, we define $\rb_i(X,Y)$. For the intuition of the construction, note that faces of $I_G(X\vvert Y)$ are order isomorphic with independent sets of $V(G)$ containing $X$ and disjoint with $Y$.\par
Suppose $v$ is a vertex of $G$, take $K=I(G)$, $K'=I(G-v)$ and $K''=I(G-N(v))$ in (\ref{recursion}). Then $K=K'\cup K''$ and $L=I_G(v\vvert\emptyset)$, so we have
\[\rb_i(G)=\rb_i(G(\emptyset\vvert v))+\rb_i(G-N(v))-\rb_i(G(v\vvert\emptyset))+\beta(N_i)+\beta(N_{i-1}).\]
Note that $I(G-N(v))$ is the collection of all simplices of the form $va_0\cdots a_p$ where $a_0\cdots a_p$ is a simplex of $I_G(v\vvert\emptyset)$, along with all faces of such simplices. That is, $I(G-N(v))$ is the cone on $I_G(v\vvert\emptyset)$ with vertex $v$, or denoted by $v*I_G(v\vvert\emptyset)$. It is an elementary fact in topology that a cone has zero reduced homology groups\cite{Mun84}:
\[\tilde{H}_i(v*I_G(v\vvert\emptyset))=0\qquad for\ all\ i.\]

That is, we have the following proposition:
\begin{prop}\label{prop:iso}
If $H$ has an isolated vertex, then $b(H)=0$.
\end{prop}

So the above equation is reduced to
\begin{equation}\label{eqn:recur1}
\rb_i(G)=\rb_i(\emptyset\vvert v)-\rb_i(v\vvert\emptyset)+\beta(N_i)+\beta(N_{i-1}), \quad\forall i. \end{equation}

Similarly, if we replace $G$ by $G-N[X]-Y$ for any vertex set $X$ and $Y$, we have 
\begin{equation}\label{eqn:recur2}
\rb_i(X,Y)=\rb_i(X\vvert Y\cup\{v\})
-\rb_i(X\cup\{v\}\vvert Y)+\beta(N'_i)+\beta(N'_{i-1}), \quad \forall i
\end{equation}
Here $N'_i$ is a subgroup of $\rH_i(I_G(\Xv v))$. We have 
\begin{equation}\label{eqn:Ni}
 \beta(N_i)\le \rb_i(v\vvert\emptyset) \mbox{ and } \beta(N'_i)\le \rb_i(\Xv v).
\end{equation}

Our proof of the main result are based on the above recursion formulas.

\section{Proof of Main Theorem}

We are going to prove \textbf{Theorem~\ref{thm:main}}: If $b(G)\ge 2$ and $b(H)\le 1$ for every induced subgraph $H$ of $G$, then $G=C_{3k}$ for some integer $k$.

In the remaining part of the paper, we fix $G$ to the one in the main result, unless it is specified individually.

Suppose $G(X|Y)$ is well defined and not a null graph. By assumption, we have $b(X\vvert Y)$ being $0$ or $1$ if $X\cup Y\not=\emptyset$. If $b(H)=1$ for some graph $H$, we denote $d(H)$ to be the dimension of the reduced Betti number taking value $1$. That is $d(G)=i$ if $\rb_i(G)=1$ and $\rb_j(G)=0$ for $j\neq i$. If $b_G=0$ then we denote $d(G)$ to be '$*$'. Note that for $b_G\geq2$, $d(G)$ is not defined. For simplicity, we write $d(X\vvert Y)$ for $d(G(X\vvert Y))$. Specially, when $X$ is not independent, we let $d(X\vvert Y)=0$, and if $G(X|Y)$ is a null graph, then $d(X\vvert Y)=-1$. That is, we consider the empty simplex as the $(-1)$-dimension sphere.

\begin{lem}\label{lem:tri-relation}
	For any disjoint vertex set $X$ and $Y$ in $G$ with $X\cup Y\not=\emptyset$ and a vertex $v$ not in $X$ or $Y$, the triple $(d(X\vvert Y),d(\Xv v),d(\vY v))$ fits into one of the following four patterns: $(k,*,k)$, $(*,*,*)$, $(*,k,k)$ and $(k+1,k,*)$ for some integer $k$.
	
	Or we can put it in triangle diagram as below:
	
\begin{center}
\begin{tabular}{ c || c  : c : c : c}
$(X,Y)$ & $k$ & $*$ & $*$ & $k+1$\\
& & & &\\
$(X\cup\{v\}\vvert Y)$\quad $(X\vvert Y\cup\{v\})$ \qquad& $*$ \qquad $k$  &  $*$\qquad $*$  &  $k$\qquad $k$  &  $k$\qquad $*$
\end{tabular}
\end{center}
\end{lem}

As the picture shown above, the triangle on the left takes one of the four patterns on the right. We say such triples(triangles) are \textit{legal} and others are \textit{illegal}. Note that if we know two corners of a legal triangle then we can determine the third.

\begin{proof}[Proof of Lemma~\ref{lem:tri-relation}]
    When $X\cup\{v\}$ is not independent, then $d(X\cup\{v\}\vvert Y)=*$, and $G-N[X]-\{v\}-Y=G-N[X]-Y$ hence $d(X\vvert Y\cup \{v\})=d(X\vvert Y)$. The triple must be one of $(k,*,k)$ or $(*,*,*)$. So we may assume $X\cup\{v\}$ is an independent set, and $v$ is a vertex in $G[X\vvert Y]$.
    
    If $V(G(X\vvert Y))=\{v\}$, and the triple is $(*,-1,-1)$. If $V(G(X\vvert Y))\not=\{v\}$, and $G(X\cup\{v\}\vvert Y)$ is null, then $v$ adjacent to all vertices in $G(X\vvert Y)$. By the minimality of $G$, we may assume $G$ does not contain triangles, hence $G(X\vvert Y)$ is a star center at $v$. As $b_0(X\vvert Y)$ is the number of components of $I(X\vvert Y)$ and $v$ is isolated in $I(X\vvert Y)$ and $b(X\vvert Y)=0$. Furthermore, $G(X\vvert Y\cup\{v\})$ is a graph with no edge, and $b(X\vvert Y\cup \{v\})=*$ by Proposition~\ref{prop:iso}. Therefore the triple will be $(0,-1,*)$.So we may assume $G(X\cup\{v\}\vvert Y)$ and $G(X\vvert Y\cup\{v\})$ are not null. 
    
	By (\ref{eqn:Ni}), if $d(\Xv v)=*$, then $\beta(N'_i)=0$ for all $i$ , and by (\ref{eqn:recur2}), $\rb_i(X\vvert Y)=\rb_i(G(\vY v))$; Thus $(d(X\vvert Y),d(\Xv v),d(\vY v))=(k,k,*)$ for some $k$ or $(*,*,*)$. 
	
	If $d(\Xv v)=k$ then $\beta(N'_k)=0$ or $1$ and $\beta(N'_i)=0$ for $i\neq k$. In the case $\beta(N'_k)=0$, by (\ref{eqn:recur2}), $\rb_k(X\vvert Y)=\rb_k(\vY v))-1$, which should be non-negative. So $\rb_k(G(\vY v))=1$ and $\rb_k(X\vvert Y)=0$, and also $\rb_i(X\vvert Y)=\rb_i(\vY v))-\rb_i(\Xv v)=0$ for $i\neq k$, hence $(d(X,Y),d(\Xv v),d(\vY v))=(*,k,k)$.
	
	In the case that $\beta(N'_k)=1$, by (\ref{eqn:recur2}) we have 
	\begin{align*}
	\rb_{k+1}(X\vvert Y)&=\rb_{k+1}(\vY v)-\rb_{k+1}(G(\Xv v))+\beta(N'_{k+1})+\beta(N'_k)\\
	&=\rb_{k+1}(\vY v)+1
	\end{align*}
	So $\rb_{k+1}(X\vvert Y)=1$, and $\rb_i(X\vvert Y)=0$ for $i\neq k+1$. Now that we have known $d(X,Y)$, $d(\Xv v)$ and all $\beta(N'_i)$'s, from (\ref{eqn:recur2}) we have $\rb_i(G(\vY v))=0$ for all $i$. So $(d(X\vvert Y),d(\Xv v),d(\vY v))=(k+1,k,*)$.
\end{proof}

\begin{lem}\label{lem:extend}
	Suppose $X, Y$ are vertex set of $G$ with $d(X\vvert Y)=k$ for some integer $k$. If $v_1, v_2$ are two vertices not in $X\cup Y$ with $d(\Xv {v_1}\}=k-1$ and $d(\Xv {v_2})=*$, then $d(\Xv {v_1,v_2})=*$.
\end{lem}

\begin{proof}
    By Lemma~\ref{lem:tri-relation}, we should have 
    $$(d(X\vvert Y),d(\Xv {v_1}),d(\vY{v_1})=(k, k-1, *)$$ and
    $$(d(X\vvert Y),d(\Xv {v_2}),d(\vY{v_2})=(k, *, k).$$
    
    Suppose $d(\Xv {v_1,v_2})\not=*$, then 
    $$(d(\Xv {v_1}),d(\Xv{v_1,v_2}),d(\XvY {v_1}{v_2})=(k-1,k-2,*)$$ and  $$(d(\Xv {v_2}),d(\Xv{v_1,v_2}),d(\XvY {v_2}{v_1})=(*, k-2, k-2).$$
	
	Now to calculate $d(\vY {v_1,v_2})$, we should have $$(d(\vY {v_1}), d(\XvY {v_2}{v_1}, d(\vY{v_1,v_2})=(*,k-2, k-2)$$ and 
	$$(d(\vY {v_2}), d(\XvY {v_1}{v_2}), d(\vY{v_1,v_2})=(k,*,k),$$ which conflict at the value of $d(\vY {v_1,v_2})$.
\end{proof}

\begin{clm}\label{clm:dG}
There is some $k\geq0$ such that $\rb_k(G)=2$ and $\rb_i(G)=0$ for all $i\neq k$. Furthermore, for every vertex $v$, $d(v\vvert\emptyset)=k-1$ and $d(\emptyset\vvert v)=k$.
\end{clm}

\begin{proof}

    Note that $N_{-1}=0$ in (\ref{eqn:recur1}), so for any vertex $v$,
	\begin{align*}
	b(G)&=b(\emptyset\vvert v)-b(v\vvert\emptyset)+2\Sigma_{i\geq0}\beta(N_i)\\
	&\leq b(\emptyset\vvert v)-b(v\vvert\emptyset)+2b(v\vvert \emptyset)\\
	&=b(\emptyset\vvert v)+b(v\vvert\emptyset).
	\end{align*}
	By the assumption of $G$, we have have $b(\emptyset\vvert v)=b(v\vvert\emptyset)=1$ and $b(G)=2$. Also, we must have $\beta(N_i)=\rb_i(G(v\vvert\emptyset))$ for all $i$ and (\ref{eqn:recur1}) is reduced to
	\begin{equation}\label{eqn:recur3}
	\begin{split}
	\rb_i(G)&=\rb_i(G(\emptyset\vvert v))-\rb_i(G(v\vvert\emptyset))+\beta(N_i)+\beta(N_{i-1})\\
	&=\rb_i(G(\emptyset\vvert v))+\rb_{i-1}(G(v\vvert\emptyset)).
	\end{split}
	\end{equation}
	
	Suppose $\rb_k(G)=\rb_l(G)=1$ for some integers $k, l$ with $k<l$. By (\ref{eqn:recur3}), for each vertex $v$, either 
    	$$v\in V_1=\{u:d(u\vvert \emptyset)=k-1, d(\emptyset\vvert u)=l\},$$ 
	or 
	    $$v\in V_2=\{u:d(u\vvert \emptyset)=l-1, d(\emptyset\vvert u)=k\}.$$ 
	
	We claim that for disjoint subsets $X$, $Y$ of $V_1$ such that $X\cup Y\neq\emptyset$, we have
	\[d(X\vvert Y)=
	\begin{cases}
	k-\vert X\vert,&Y=\emptyset,\\
	*,&X,Y\neq\emptyset,\\
	l,&Y=\emptyset.
	\end{cases}
	\]
	Or we can put it in triangle diagram as below:
	\begin{center}
		\begin{tikzpicture}
		\node at (-4-9+1/4,-1/4) {$(\vert X\vert,\vert Y\vert)$:};
		\foreach \n in {1,...,3} {
			\foreach \k in {0,...,\n} {
				{
					\pgfmathsetmacro\tuple{\n-\k}
					\node at (\k-\n/2-9+1/4,-\n/4*3) {$(\pgfmathprintnumber\tuple,\k)$};
				}
			}
		}
		\node at (0-9+1/4,-3) {$\vdots$};
		\node at (-5/2-9+1/4,-5/4*3) {$(t,0)$};
		\node at (-3/2-9+1/4,-5/4*3) {$(\!t\!\!-\!\!1\!,\!1\!)$};
		\node at (-1/2-9+1/4,-5/4*3) {$\cdots$};
		\node at (1/2-9+1/4,-5/4*3) {$\cdots$};
		\node at (3/2-9+1/4,-5/4*3) {$(\!1\!,\!t\!\!-\!\!1\!)$};
		\node at (5/2-9+1/4,-5/4*3) {$(0,t)$};
		\draw (-11/2,0) -- (-11/2,-4);
		\node at (-4-1/4,-1/4) {$d(X\vvert Y)$:};
		\foreach \n in {1,2,3} {
			\node at (-\n/2-1/4,-\n/4*3) {$k\!-\!\n$};
			\node at (\n/2-1/4,-\n/4*3) {$l$};
		}
		\node at (0-1/4,-3/2) {$*$};
		\node at (-1/2-1/4,-9/4) {$*$};
		\node at (1/2-1/4,-9/4) {$*$};
		\node at (0-1/4,-3) {$\vdots$};
		\node at (-5/2-1/4,-5/4*3) {$k\!-\!t$};
		\node at (-3/2-1/4,-5/4*3) {$*$};
		\node at (-1/2-1/4,-5/4*3) {$\cdots$};
		\node at (1/2-1/4,-5/4*3) {$\cdots$};
		\node at (3/2-1/4,-5/4*3) {$*$};
		\node at (5/2-1/4,-5/4*3) {$l$};
		\end{tikzpicture}
	\end{center}
	
	We prove by induction on $t=|X\cup Y|$. It is true for $t=1$ by definition of $V_1$.

	Suppose we have proved for $\vert X\cup Y\vert\leq t-1$ for $t\geq 2$. Note that for the first $(t-1)$ rows, $d(X\vvert Y)$ is determined by $(|X|,|Y|)$, so we may use $d(|X|,|Y|)$ to denote $d(X\vvert Y)$. In $t$-row with $W=X\cup Y$ given, by repeatedly using Lemma~\ref{lem:tri-relation}, the $t$-th row are determined by $d(W\vvert \emptyset)$, and $d(X\vvert Y)$ are also determined by $(|X|,|Y|)$. Furthermore, the triple $(d(X\vvert Y),d(X\cup\{v\}\vvert Y),d(X\vvert Y\cup\{v\}))$ in Lemma~\ref{lem:tri-relation} can also be replaced by $(d(|X|,|Y|),d(|X|+1,|Y|), d(|X|,|Y|+1))$, which form a small triangle in the triangle diagram above. 
	
	By Lemma~\ref{lem:tri-relation}, there are at most two possible lists of the values on the $t$-th row, depending on $d(t,0)$. They are $(k-t, *, \dots, *, l)$ or $(*, k-t+1, k-t+1, \dots, k-t+1, *)$. The later one is legal only when $k=l+t-2$, which conflict with the assumption that $k<l$. Therefore the $t$-th row must be $(k-t, *,\cdots, *, l)$ for any $W=X\cup Y$ with size $t$. By induction, the claim is true for all rows. 
	
	Specially, we have $d(\emptyset\vvert V_1)=l$. Using the same argument above with $t=2$, as $l\neq k+t-2=k$, we can get for $\forall u, v\in V_2$, we have $d(u,v\vvert\emptyset)=l-2$, which implies that any two vertices $u, v$ in $V_2$ are not adjacent, hence $V_2$ is an independent set. However, as $G(\emptyset \vvert V_1)=G[V_2]$, by Proposition~\ref{prop:iso}, $b(\emptyset \vvert V_1)=0$, contradicting that $d(\emptyset\vvert V_1)=l$.
	
	So we have $\rb_k(G)=2$ for some $k$. By (\ref{eqn:recur3}), we have  $d(v\vvert\emptyset)=k-1$ and $d(\emptyset\vvert v)=k$ for every vertex $v$. 
\end{proof}

Throughout the rest of this article, we use the letter $k$ to refer to the integer we obtained in the above theorem.

Let $u$, $v$ be two vertices of $G$. As $d(u\vvert\emptyset)=k-1$, by Lemma~\ref{lem:tri-relation}, $(d(u\vvert\emptyset),d(u,v\vvert\emptyset),d(u\vvert v))$ is either $(k-1,k-2,*)$ or $(k-1,*,k-1)$. We construct a new graph $H$ on $V(G)$ such that $u$, $v$ are adjacent if and only if $d(u,v\vvert\emptyset)=k-2$. The following propositions of $H$ immediately follow:

\begin{prop} \label{prop:xyinH}
In $H$, any two vertices $u, v$ satisfies
\begin{enumerate}
    \item If $u\sim v$ in $H$, then $u\nsim v$ in $G$. That is, $E(G)\cap E(H)=\emptyset$. 
    \item If $u\sim v$ in $H$, then $d(u,v\vvert \emptyset)=k-2$, $d(u\vvert v)=d(v\vvert u)=*$, and $d(\emptyset\vvert u,v)=k$, 
    \item If $u\nsim v$ in $H$, then $d(u,v\vvert \emptyset)=d(\emptyset\vvert u,v)=*$, and $d(u\vvert v)=d(v\vvert u)=k-1$.
\end{enumerate}
\end{prop}

The following proposition is a key feature of $H$.

\begin{lem}\label{1comp}
	Every component $C$ of $H$ is a complete graph. Furthermore, for any disjoint subsets $X$ and $Y$ of $V(C)$ with $X\cup Y\neq\emptyset$, we have
	\[d(X\vvert Y)=
	\begin{cases}
	k-\vert X\vert,&Y=\emptyset,\\
	*,&X,Y\neq\emptyset,\\
	k,&X=\emptyset.
	\end{cases}
	\]
	Or we can put it in triangle diagram as below:
	\begin{center}
		\begin{tikzpicture}
		\node at (-4-9+1/4,-1/4) {$(\vert X\vert,\vert Y\vert)$:};
		\foreach \n in {1,...,3} {
			\foreach \k in {0,...,\n} {
				{
					\pgfmathsetmacro\tuple{\n-\k}
					\node at (\k-\n/2-9+1/4,-\n/4*3) {$(\pgfmathprintnumber\tuple,\k)$};
				}
			}
		}
		\node at (0-9+1/4,-3) {$\vdots$};
		\node at (-5/2-9+1/4,-5/4*3) {$(t,0)$};
		\node at (-3/2-9+1/4,-5/4*3) {$(\!t\!\!-\!\!1\!,\!1\!)$};
		\node at (-1/2-9+1/4,-5/4*3) {$\cdots$};
		\node at (1/2-9+1/4,-5/4*3) {$\cdots$};
		\node at (3/2-9+1/4,-5/4*3) {$(\!1\!,\!t\!\!-\!\!1\!)$};
		\node at (5/2-9+1/4,-5/4*3) {$(0,t)$};
		\draw (-11/2,0) -- (-11/2,-4);
		\node at (-4-1/4,-1/4) {$d(X,Y)$:};
		\foreach \n in {1,2,3} {
			\node at (-\n/2-1/4,-\n/4*3) {$k\!-\!\n$};
			\node at (\n/2-1/4,-\n/4*3) {$k$};
		}
		\node at (0-1/4,-3/2) {$*$};
		\node at (-1/2-1/4,-9/4) {$*$};
		\node at (1/2-1/4,-9/4) {$*$};
		\node at (0-1/4,-3) {$\vdots$};
		\node at (-5/2-1/4,-5/4*3) {$k\!-\!t$};
		\node at (-3/2-1/4,-5/4*3) {$*$};
		\node at (-1/2-1/4,-5/4*3) {$\cdots$};
		\node at (1/2-1/4,-5/4*3) {$\cdots$};
		\node at (3/2-1/4,-5/4*3) {$*$};
		\node at (5/2-1/4,-5/4*3) {$k$};
		\end{tikzpicture}
	\end{center}
\end{lem}

\begin{proof}
	Suppose $C$ is not complete, there must exist three distinct vertices in $C$ $u$, $v$ and $w$ such that $u\sim v$, $v\sim w$ but $u\nsim w$ in $H$. Since $d(u\vvert\emptyset)=k-1, d(u,v\vvert \emptyset)=k-2, d(u,w\vvert \emptyset)=*$, by Lemma~\ref{lem:extend}, $d(u,v,w\vvert\emptyset)=*$. So the triple $(d(u,v\vvert\emptyset),d(u,v,w\vvert\emptyset),d(u,v\vvert w))=(k-2,*,k-2)$.  But as $d(\emptyset\vvert w)=k, d(u\vvert w)=k-1,d(v\vvert w))=*$, by Lemma~\ref{lem:extend}, we should have $d(u,v\vvert w)=*$. Contradiction! 
	
	So $C$ must be complete, and it will imply the first two rows of the triangle diagram. The remaining level can be proved inductively just like to the one in Theorem~\ref{clm:dG}, with $k=l$ and $k\neq l+t-2$ when $t\ge 3$.
\end{proof}

The following result follows immediately.

\begin{clm}\label{clm:ngb}
	There does not exist a vertex $v$ with all neighbors in $G$ located in one component of $H$.
\end{clm}

\begin{proof}
	Suppose there is a component $C$ of $H$ such that $N_G(v)\subseteq C$. By Lemma~\ref{1comp}, $d(\emptyset\vvert C)=d(G-C)=k$, but $v$ is an isolated vertex in $G-C$, we have $b(G-C)=0$ by Proposition~\ref{prop:iso}. Contradiction. 
	
\end{proof}

\begin{lem}\label{lem:2edge4comp}
	There do not exist two edges $v_1v_2,v_3v_4$ in $G$, with $v_1,v_2,v_3,v_4$ located in four distinct components of $H$.
\end{lem}

\begin{proof}
	 Suppose they are located in distinct components of $H$. We consider $d(X\vvert Y)$ for disjoint subsets $X, Y$ in $\{v_1,v_2,v_3,v_4\}$ with $X\cup Y\not=\emptyset$. We claim to have the following triangle diagram:
	\begin{center}
		\begin{tikzpicture}
		\node at (-4-4,-1/4) {$(\vert X\vert,\vert Y\vert)$:};
		\foreach \n in {1,...,4} {
			\foreach \k in {0,...,\n} {
				{
					\pgfmathsetmacro\tuple{\n-\k}
					\node at (\k-\n/2-4,-\n/4*3) {$(\pgfmathprintnumber\tuple,\k)$};
				}
			}
		}
		\draw (-5/4,0) -- (-5/4,-13/4);
		\node at (-4+4,-1/4) {$d(X\vvert Y)$:};
		\node at (-1/2+4,-3/4) {$k\!-\!1$};
		\node at (1/2+4,-3/4) {$k$};
		\node at (-1+4,-3/2) {$*$};
		\node at (0+4,-3/2) {$k\!-\!1$};
		\node at (1+4,-3/2) {$*$};
		\node at (-3/2+4,-9/4) {$*$};
		\node at (-1/2+4,-9/4) {$*$};
		\node at (1/2+4,-9/4) {$k\!-\!1$};
		\node at (3/2+4,-9/4) {$k\!-\!1$};
		\node at (-2+4,-3) {$*$};
		\node at (-1+4,-3) {$*$};
		\node at (0+4,-3) {$*$};
		\node at (1+4,-3) {$?$};
		\end{tikzpicture}
	\end{center}
	The first row is implied by Claim~\ref{clm:dG}. The second row follows from the assumption that each pair of the $v_i$'s belongs to different components of $H$. Note that when $|X|\ge 3$, $X$ is not an independent set, hence $d(X\vvert Y)=*$. Therefore the first term of the third row and fourth row are $*$. And we can get the rest of the third row using Theorem~\ref{lem:tri-relation}. Similarly the first three terms of the fourth row are '$*$'. But there is no proper value for $d(X\vvert Y)$ with $(|X|,|Y|)=(1,3)$ that fits Lemma~\ref{lem:tri-relation}.
\end{proof}

Now we are ready to complete the proof of the main theorem.

\begin{proof}[Proof of Theorem~\ref{thm:main}]
	As $b(C_{3k})\ge |f_{C_{3k}}|=2$ for any $k$, we just need to show that $G$ contains an induced $C_{3k}$ for some $k$.
	
	First each component of $H$ is an independent set in $G$. By Lemma~\ref{clm:ngb}, the neighbors of any vertex in $G$ are located in at least two components of $H$; and by Lemma~\ref{lem:2edge4comp}, there is no two edges with all ends in four distinct components. 
	With these together it is easy to deduce that $H$ have have exactly three components $C_0$, $C_1$, $C_2$, and every vertex in $C_i$ has neighbours in $C_{i-1}$ and $C_{i+1}$ for each $i$ (indices modulo 3). We orient all the edges of $G$ from $C_i$ to $C_{i+1}$ for $i=0, 1, 2$. Then every vertex has positive out-degree, therefore there is an induced directed cycle $[v_1,v_2,\cdots, v_r]$, which must has length divisible by three. 
\end{proof}

\section{Acknowledgement}

The first author will like to thank his colleagues Guozhen Wang and Xiping Zhang in Shanghai Center for Mathematical Sciences for sharing knowledge in Algebraic Topology, and to thank Qiqin Xie for the early discussion on this problem.

\end{document}